\documentclass[12pt]{amsart}
\font\emailfont=cmtt10

\usepackage{amsmath,amsthm,amsfonts,amscd,flafter,epsf}

\newcommand\commentable[1]{#1}

\newcommand\Id{\mathrm{Id}}

\newtheorem{theorem}{Theorem}[section]
\newtheorem{prop}[theorem]{Proposition}

\newtheorem{lemma}[theorem]{Lemma}

\newtheorem{defn}[theorem]{Definition}

\newtheorem{remark}[theorem]{Remark}

\def\endproof{\relax\ifmmode\expandafter\endproofmath\else
  \unskip\nobreak\hfil\penalty50\hskip.75em\hbox{}\nobreak\hfil\bull
  {\parfillskip=0pt \finalhyphendemerits=0 \bigbreak}\fi}
\def\endproofmath$${\eqno\bull$$\bigbreak}
\def\bull{\vbox{\hrule\hbox{\vrule\kern3pt\vbox{\kern6pt}\kern3pt\vrule}\hrule}}
\newcommand{\smargin}[1]{}

\newcounter{bean}
\newcommand{\Q}{\mathbb{Q}}
\newcommand{\R}{\mathbb{R}}

\newcommand{\Z}{\mathbb{Z}}

\newcommand{\Zmod}[1]{\Z/{#1}\Z}

\newcommand{\Ker}{\mathrm{Ker}}

\newcommand{\cm}{\cdot}

\newcommand{\CDisk}{D}

\newcommand{\ModSWfour}{\mathcal{M}}
\newcommand{\ModFlow}{\ModSWfour}

\newcommand{\SpinC}{{\mathrm{Spin}}^c}

\newcommand\abuts\Rightarrow
\newcommand\Sym{\mathrm{Sym}}

\newcommand\HFpRed{\HFp_{\red}}

\newcommand\RelSpinC{\underline{\SpinC}}
\newcommand\relspinc{\underline{\spinc}}

\newcommand\Filt{\mathcal F}

\newcommand\x{\mathbf x}

\newcommand\y{\mathbf y}

\newcommand\ModSphere{\ModFlow\left({\mathbb S}\longrightarrow 
\Sym^{g-1}(\Sigma_{1})\times \Sym^2(\Sigma_{2})\right)}
\newcommand\ModSpheres\ModSphere

\newcommand\CFa{\widehat{CF}}

\newcommand{\red}{\mathrm{red}}

\newcommand\HFp{\HFb}

\newcommand\HFinf{HF^\infty}

\newcommand\HFa{\widehat{HF}}
\newcommand\HFb{HF^+}

\newcommand\Mas{\mu}
\newcommand\UnparModSp{\widehat \ModSp}
\newcommand\UnparModFlow\UnparModSp
\newcommand\Mod\ModSp

\newcommand{\cald}{{\mathcal D}}

\newcommand{\spinc}{\mathfrak s}

\newcommand\ModMaps{\mathcal M}
\newcommand\ModSp\ModMaps

\newcommand\Ta{{\mathbb T}_{\alpha}}
\newcommand\Tb{{\mathbb T}_{\beta}}
\newcommand\Tc{{\mathbb T}_{\gamma}}

\newcommand\alphas{\mbox{\boldmath$\alpha$}}

\newcommand\betas{\mbox{\boldmath$\beta$}}
\newcommand\gammas{\mbox{\boldmath$\gamma$}}

\newcommand\uHFa{\underline{\HFa}}

\newcommand\Dual{\mathcal D}
\newcommand\Duality\Dual

\newcommand\RightDehn{R}
\newcommand\relspinct{\underline{\mathfrak t}}
\newcommand\Mass{\overline n}
\newcommand\spinccan{{\mathfrak k}}
\newcommand\spinccant{\ell}

\newcommand\FiltZero{\mathcal K}

\newcommand\MClass{\mathcal M}

\newcommand\CFK{CFK}

\newcommand\CFKa{\widehat\CFK}

\newcommand\CFKinf{\CFK^{\infty}}

\newcommand\Mark{m}

\commentable{

\title[{Heegaard Floer homologies and contact structures}]
{Heegaard Floer homologies and contact structures}

\author[Peter Ozsv{\'a}th]{Peter Ozsv\'ath}
\address{Department of
Mathematics, Columbia University, New York 10027 \newline
\indent{\emailfont{petero@math.columbia.edu}}}

\author[Zolt{\'a}n Szab{\'o}]{Zolt{\'a}n Szab{\'o}} 
\address{Department of
Mathematics, Princeton University, New Jersey 08540 \newline
\indent{\emailfont{szabo@math.princeton.edu}}}}

\begin{document}
\begin{abstract}
Given a contact structure on a closed, oriented three-manifold $Y$, we
describe an invariant which takes values in the three-manifold's Floer
homology $\HFa$ (in the sense of~\cite{HolDisk}).  This invariant
vanishes for overtwisted contact structures and is non-zero for Stein
fillable ones.  The construction uses of Giroux's interpretation of
contact structures in terms of open book decompositions
(see~\cite{Giroux}), and the knot Floer homologies introduced
in~\cite{Knots}. 
\end{abstract}

\maketitle
\section{Introduction}

In~\cite{HolDisk}, we defined several Floer homology groups associated
to a closed oriented three-manifold $Y$. The simplest variant of these
groups is $\HFa(Y)$ (which can be further decomposed according to
$\SpinC$ structures over $Y$).  Our goal here is to associate to each
contact structure over $Y$ a corresponding element in this group
(actually, the element strictly speaking lies in $\HFa(-Y)$).  The
motivation to find such an invariant came from an analogous
construction due to Kronheimer and Mrowka (see~\cite{KMcontact}),
which can be interpreted as giving a map from the set of contact
structures over $Y$ to the ``Seiberg-Witten Floer cohomology'' of
$Y$. Although the end results of these two constructions are very
similar, their specifics are quite different. While Kronheimer and
Mrowka use the contact structure to give a suitable boundary condition
for the Seiberg-Witten monopoles over $[0,\infty)\times Y$, our
constructions here proceed by analyzing the open book decomposition
induced from the contact structure, as provided by recent work of
Giroux~\cite{Giroux} (see also,~\cite{Torisu}, \cite{GirouxPaper}).
Before proceeding to the details, we recall some background.

Informally, a contact structure $\xi$ on an oriented three-manifold
$Y$ is an everywhere totally non-integrable two-plane field. More
precisely, when $Y$ is oriented, the two-plane field $\xi$ represents a
contact structure if there is a one-form $\lambda$ so that $\xi=\Ker
\lambda$, and $\lambda\wedge d\lambda$ represents a volume form for
$Y$ with the specified orientation. (These are often called {\em
cooriented} contact structures in the literature; since they are the
only contact structures we consider in this paper, we drop the
modifier from our terminology.) Two nowhere integrable two-plane
fields $\xi_1$ and $\xi_2$ represent the same contact structure if
there is a diffeomorphism $\phi$ of $Y$ which is isotopic to the
identity map and which carries $\xi_1$ to $\xi_2$.  (According to a
basic theorem of Gray~\cite{Gray}, this is equivalent to the condition that
the two-plane fields $\xi_1$ and $\xi_2$ are homotopic through nowhere
integrable two-plane fields.)

Informally, an open book decomposition of $Y$ is a non-empty,
oriented link $L\subset Y$ (the {\em binding}) whose complement fibers over the
circle. More precisely, let $\phi$ be an automorphism of an oriented
surface $F$ with non-empty boundary, and suppose that $\phi$ fixes
$\partial F$.  We can form the mapping torus $M_\phi=F\times
[0,1]/(x,0)\sim (\phi(x),1)$ to obtain a three-manifold which
naturally fibers over the circle, and whose boundary is $\partial
F\times S^1$.  There is a canonically associated closed three-manifold
$Y_0$
obtained from $M_\phi$ by attaching solid tori $\partial F\times
\CDisk$ using the identifications suggested by the notation. The data
$(F,\phi)$ is called an open book decomposition of $Y$, and $\phi$ is
called the monodromy of the open book. The binding here corresponds to
the boundary components of $F$.  When the link has a single component,
the open book decomposition is called a fibered knot.  A construction
of Thurston and Winkelnkemper (see~\cite{ThurstonWinkelnkemper})
associates to an open book decomposition of $Y$ a contact
structure. Indeed, according to recent work of Giroux ~\cite{Giroux},
every contact structure is induced from an open book in this way and,
in fact, Giroux gives an explicit criterion for when two open books
induce the same contact structure.  (see also
Subsection~\ref{subsec:Giroux} below for a more detailed discussion of
these moves)

To describe the Floer homology class associated to $\xi$, we use the
knot invariants introduced in~\cite{Knots} (compare
also~\cite{Rasmussen}).  When a three-manifold is equipped with an
oriented knot $K$ and a Seifert surface $F$ for $K$, there is an
induced $\Z$-filtration on its Floer complex $\CFa$ which calculates
$\HFa$. We let $\FiltZero(Y,K,F,m)\subset \CFa(Y)$ denote the
subcomplex generated by all intersection points with filtration level
$\leq m$.  We review the construction of the knot filtration in
Section~\ref{sec:TopPre}, including a review of some of the necessary
topological preliminaries. We will typically drop $F$ from the notation
when it is clear from the context (for instance, when 
the knot arises from an open book decomposition).

To conform with the usual conventions in contact geometry, we will
typically work with Floer cohomology rather than homology or,
equivalently, the Floer homology of the three-manifold $-Y$ whose
orientation is the opposite of the one given by the contact structure.

\begin{theorem}
\label{intro:FiberedKnots}
Let $Y$ be an oriented three-manifold, $K\subset Y$ be a fibered knot
of genus $g$.  Then, $H_*(\FiltZero(-Y,K,F,-g))\cong \Z$.
\end{theorem}

The above theorem allows us to define the Floer homology class
associated to a fibered knot $K\subset Y$.

\begin{defn}
Let $K\subset Y$ be a fibered knot, and let 
$$c_0(K)\in
H_*(\FiltZero(-Y,K,-g))\cong \Z$$ be a generator.  Then, the
invariant of the fibered knot, $c(K)\in \HFa(-Y)$, is the image of
$c_0(K)$ under the map on homology induced from the natural inclusion
map
$$\iota\colon \FiltZero(-Y,K,-g)\longrightarrow
\CFa(-Y).$$
\end{defn}

In the above definition, $c(K)$ is well-defined only up to sign. Of
course, if we perform the constructions here with $\Zmod{2}$
coefficients, then this sign ambiguity disappears. Otherwise, we are
to view $c(K)$ as an element of the set $\HFa(-Y)/(\pm 1)$.

In the course of proving Theorem~\ref{intro:FiberedKnots}, we also
establish naturality properties of the isomorphism. These properties,
together with Giroux's result, lead to the following result
established in Section~\ref{sec:Equivalence}:

\begin{theorem}
\label{intro:SameContact}
Suppose $K_1$ and $K_2$ are a pair of knots in $Y$ which represent the
same contact structure $\xi$, then, the invariants $c(K_1)$
and $c(K_2)$ coincide. 
\end{theorem}

In view of the above result, if $\xi$ is a contact structure over $Y$,
we write $c(\xi)$ for the element defined earlier. In
Subsection~\ref{sec:Equivalence}, we also discuss how to extract
classical data about the homotopy class of the two-plane field from
the Floer homology class $c(\xi)$. In particular, a two-plane field
$\xi$ induces a $\SpinC$ structure $\spinc(\xi)$ on $Y$, and we show that the
element $c(\xi)$ is supported in the corresponding summand
$\HFa(Y,\spinc(\xi))\subset \HFa(Y)$.

Recall that a contact structure is called {\em overtwisted} if there
is a disk in $Y$ which is transverse to $\xi$ in a neighborhood of the
boundary, but whose boundary is tangent to $\xi$
(see~\cite{EliashbergOvertwisted}). A fundamental result of Eliashberg
states that each homotopy class of two-plane field over $Y$ contains
exactly one isotopy class of overtwisted contact structure.  A contact
structure which is not overtwisted is called tight.

\begin{theorem}
\label{intro:OverTwisted}
Suppose that $\xi$ is overtwisted. Then, the induced element
$c(\xi)\in\HFa(-Y)/(\pm 1)$ is trivial.
\end{theorem}

Another theorem of Eliashberg and Gromov (see~\cite{EliashbergFilling}
and~\cite{Gromov}, see also~\cite{KMcontact} for a Seiberg-Witten
proof) states that symplectically semi-fillable contact structures are
tight (see those references for the definition).  Although presently
the symplectic semi-fillability hypothesis cannot be described
explicitly in terms of the monodromy map, there is a special class of
such contact structures which, thanks again to the work of
Giroux~\cite{Giroux}, admits such a concrete characterization: the
Stein fillable ones.

Recall (c.f.~\cite{Gromov}) that a Stein surface is a complex surface 
which admits a Morse
function $$f\colon S\longrightarrow [0,1]$$ with $\partial
S=f^{-1}(1)$ such that away from the critical points of $f$, the field
of complex tangencies to each preimage $f^{-1}(t)$ is a
positively-oriented contact structure (with respect to the orientation
induced from $\partial f^{-1}([0,t])$). The contact structure
$(Y,\xi)$ is Stein fillable if it is realized as the boundary of a
Stein surface in the above sense.  According to Giroux, these are the
contact structures whose monodromy can be expressed as a product of
right-handed Dehn twists. With this characterization, we prove the following:

\begin{theorem}
\label{intro:SteinFillable}
Suppose that $\xi$ is a Stein fillable contact structure on a closed,
oriented three-manifold $Y$. Then, the induced element
$c(\xi)\in\HFa(-Y)/(\pm 1)$ is non-trivial.
\end{theorem}

Of course, Theorems~\ref{intro:OverTwisted}
and~\ref{intro:SteinFillable} together give an alternate proof that
Stein fillable contact structures are tight.

The invariant $c(\xi)$ has an alternate description which does not
refer to the knot invariants. An open book decomposition naturally
gives rise to a cobordism $W$ from a fibered three-manifold $-Y_0$ to
the given three-manifold $-Y$, and the fibered three-manifolds with
fiber genus $g>1$ have been shown to have $\HFp(-Y_0,\spinccant)\cong
\Z$ (c.f. Theorem~\ref{HolDiskSymp:thm:ThreeManifoldsFiber}
of~\cite{HolDiskSymp}), where $\spinccant$ is the ``canonical $\SpinC$
structure'', i.e.  the one represented by a vector field which is
transverse to the fibers.  Now, if $K\subset Y$ is a fibered knot with
genus $g>1$, and ${\widehat c}(K)$ is a generator of
$\HFa(-Y_0,\spinccant)$ whose image in $\HFp(-Y_0,\spinccant)$ is a
generator, we show that $c(\xi)$ agrees with the image of ${\widehat
c}(K)$ under the map on $\HFa$ induced by the cobordism $W$.

As an application of this alternative description, we give a result
for fibered knots in the three-sphere. As motivation,
suppose that $Y$ is a three-manifold with
an open book decomposition whose binding $K$ is connected.  Clearly,
for any integer $n$, there is an induced open book decomposition of
the surgered manifold $Y_{1/n}(K)$. Indeed, it is an easy consequence
of Giroux's theory (together with the observation that the mapping
class group of a surface with a single boundary component is generated
as a monoid by right-handed Dehn twists and arbitrary Dehn twists
parallel to the boundary) that for any sufficiently large $n$, the
induced contact structure on $Y_{-1/n}(K)$ is Stein fillable.

This raises the following question: suppose that $(Y,K)$ is an open
book decomposition, what is the minimal value of $n$ for which the
open book decomposition $Y_{-1/n}(K)$ is tight? In this direction, we
have the following application of the invariants of the present paper:

\begin{theorem}
\label{intro:MinusOneFibered}
Let $K\subset S^3$ be a fibered knot in $S^3$. Then,
the contact structure on $S^3_{-1}(K)$ induced from the open book
decomposition coming from $K$ is tight.
\end{theorem}

The key property of $S^3$ used in the above result is that
$\HFpRed(S^3)=0$. Indeed, the above theorem can be generalized to the
case of any three-manifold $Y$ with $\HFpRed(Y)=0$, a class of
three-manifolds which is closed under connected sum and orientation
reversal, and includes, for example, $S^2\times S^1$, the Poincar\'e
homology sphere, and any lens space (c.f. Theorem~\ref{thm:FiberedGen} below
for the statement in this case).

\subsection{Remarks}
Given an open book decomposition for the contact structure, we exhibit
the intersection point (in an appropriate Heegaard diagram)
representing a cycle for $c(\xi)\in\HFa(-Y)$, but it is, in general,
difficult to detect whether the cycle is a boundary.  Of course, in
cases where the knot complex $\CFKinf(Y,K)$ can be explicitly
calculated, this problem is easily solved (c.f.~\cite{AltKnots} for
further examples).
However, finding an algorithmic procedure for making these calculations in
general remains a very interesting problem, both for the present
application to contact geometry, and also for three-dimensional
topology.

\subsection{Relationship with Seiberg-Witten theory}

It is natural to conjecture that the image of $c(\xi)$ in $\HFp(-Y)$
corresponds (under the conjectured identification of this latter group
with Seiberg-Witten Floer homology) to the Seiberg-Witten relative
invariant induced from the symplectization of $\xi$,
defined following~\cite{KMcontact}. In a related direction, it would be very
interesting to extend the non-vanishing result of
Theorem~\ref{intro:SteinFillable} to larger classes of contact
structures: e.g. symplectically semi-fillable ones, or perhaps even to
all tight ones.

\subsection{Acknowledgements}

The authors wish to warmly thank Ko Honda, Paolo Lisca, and Andr{\'a}s
Stipsicz for interesting discussions.

\section{Preliminaries}
\label{sec:TopPre}

\subsection{The mapping class class group and open books}

When considering open book decompositions, one considers the mapping
class of the monodromy $\phi$, i.e.  $\phi$ modulo isotopies which fix
the boundary.  Operations on mapping classes give rise to operations
on open book decompositions.

One such operation which will be particularly useful for us is the boundary
connected sum. Fix surfaces $F_1$ and $F_2$,
choose a pair of points on their boundaries $\partial F_1$ and
$\partial F_2$, and let $F_1\#_b F_2$ denote the corresponding
boundary connected sum.  If $\phi_1$ and $\phi_2$ are mapping classes
of the surfaces-with-boundary $F_1$ and $F_2$, then there is an induced
automorphism $\phi_1\#_b\phi_2$ of $F_1\#_b F_2$. Correspondingly, if
$(F_1,\phi_1)$ and $(F_2,\phi_2)$ are open book decompositions of
$Y_1$ and $Y_2$, then $(F_1\#_b F_2,\phi_1\#_b\phi_2)$ is an open book
decomposition of the connected sum $Y_1\# Y_2$.

Another useful operation is composition with Dehn twists. Recall that
the mapping class group of the annulus $A=[0,1]\times S^1$ isomorphic
to $\Z$, generated by the ``right-handed Dehn twist,'' a
diffeomorphism whose mapping class is determined by the property that
the oriented intersection number $$\#\Big([0,1]\times\{x\}\cap
\psi([0,1]\times\{x\})\Big)=-1.$$
(The other generator, of course, is
called the left-handed Dehn twist.)  More generally, fix a simple,
closed curve $\gamma$ in the interior of $F$. The {\em right-handed
  Dehn twist along $\gamma$}, $\RightDehn_\gamma$ is the mapping class
which is the identity away from an annular neighborhood of $\gamma$
and whose restriction to this neighborhood is the right-handed Dehn
twist of the annulus.  \footnote{With these conventions, the monodromy
  of a Lefschetz fibration around a singular fiber is a 
  right-handed Dehn twist.}  This has the following three-dimensional
interpretation. Suppose that $(F,\phi)$ is an open book decomposition
of $Y$, then the curve $\gamma$ inside $F$ can be viewed as a knot in
$Y$. Indeed, the tangent bundle of $F$ gives $\gamma$ a framing, and
the three-manifold $Y_{-1}(\gamma)$ obtained by surgery along $\gamma$
with framing $-1$ (with respect to the surface framing) inside $Y$
admits the open book decomposition $(F,\phi\circ\RightDehn_\gamma)$.

Consider for example the right-handed trefoil knot $T_{\ell}$ in $S^3$.
This has a fibration $(E,\phi_0)$, where $E$ is a surface with genus
one and a single boundary circle. Letting $\alpha$ and $\beta$ be a
pair of curves in $E$ which meet transversally in a single
intersection point, the monodromy $\phi_0$ is a product of right-handed
Dehn twists about $\alpha$ and $\beta$.

If $\phi$ is the identity map of a surface $F$ of genus $g$ with a
single boundary component, the corresponding open book decomposition
describes a knot $B$ inside $\#^{2g}(S^2\times S^1)$. For example,
when $g=1$, this knot is the remaining component of the Borromean link, 
after one performs $0$-surgery on the other two components.

\subsection{Giroux's stabilization theorem}
\label{subsec:Giroux}

Suppose that $Y$ has an open book decomposition $(F,\phi)$. The
mapping torus $M_\phi$ admits a two-plane field, given by the tangent
spaces to the fibers $F$. Now, there is a canonical tight contact
structure on the solid torus $D^2\times S^1$ whose 
induced characteristic foliation on $(\partial D^2)\times S^1$
consists of circles $x\times S^1$. In $Y=M\cup_{\partial F\times S^1}
\left(F\times \CDisk\right)$, these two two-plane fields agree on the overlap,
giving a continuous two-plane distribution. The Thurston-Winkelnkemper
contact structure is obtained as a $C^0$-small perturbation of this
two-plane field.

Giroux shows in~\cite{Giroux} that any contact structure over $Y$
arises from this construction. Moreover, he identifies which open
books give rise to the same contact structure, using the following
stabilization procedure. Let $F$ be a
surface-with-boundary, and choose a pair of points in its
boundary. Attach a one-handle to $F$ along the boundary to obtain a
new surface $F'$, and let $\gamma$ be any curve obtained by closing up
the core of the one-handle in the interior of $F$. We say that
$(F',\phi\circ \RightDehn_\gamma)$ is obtained from $(F,\phi)$ by a simple
Giroux stabilization. More generally, we call a Giroux stabilization a
finite sequence of simple Giroux stabilizations.  Note that Giroux
stabilizations leave the underlying three-manifold $Y$ unchanged, and
indeed they induce isomorphic contact structures over $Y$. (It is
important to underscore the fact that Giroux stabilizations involve
compositions with right-handed Dehn twists: the analogous operation
using left-handed Dehn twists {\em does} change the contact structure
on $Y$.) Giroux's result further states that, conversely, if two open
book decompositions induce the same contact structure, then they can
be connected by a sequence of Giroux stabilizations and
de-stabilizations.  This is easily seen to be equivalent to the
statement that both open books have a common Giroux stabilization.

For example, if $K$ is a fibered knot in $S^3$, then the connected sum
$K\# T_r$ is a Giroux stabilization of $K$.

For simplicity, we typically consider open book decompositions where
the binding is connected. Clearly, any open book decomposition has a
Giroux stabilization with this property.

As we mentioned earlier, a related result of Giroux (see
also~\cite{Giroux}) gives a criterion for Stein fillability of a
contact structure in terms of the monodromy of the open book.  This
result states that the induced contact structure is Stein fillable if
and only if the monodromy map $\phi$ of the open book can be expressed
as a product of right-handed Dehn twists.

\subsection{Multiply-filtered chain complexes}
\label{subsec:KnotHomology}

In~\cite{Knots}, we define an invariant for a three-manifold equipped
with an oriented, null-homologous knot $K$. We sketch the special case
of this construction which we need in the present application.

The data $(Y,K)$ gives rise to a Heegaard diagram $(\Sigma,
\alphas,\mu\cup\betas_0)$, where $\alphas$ and $\betas=\mu\cup\betas_0$ are attaching
circles for $(Y,K)$, and the curve $\mu$ is a distinguished meridian
for the knot $K$, so that $(\Sigma,\alphas,\betas_0)$ represents the
knot complement.  Indeed, we can choose a reference point $m$ to lie
on $\mu$.  Using the orientations, we obtain a doubly-pointed Heegaard
diagram $(\Sigma, \alphas,\betas,w,z)$ where $w$ and $z$ are reference
points close to the initial point $m\in\mu$, but lying on either side
of it. More precisely, if $\lambda$ is a longitude for the knot,
thought of as an oriented curve in $\Sigma$, and $\delta$ is the
oriented path from $z$ to $w$, then the algebraic intersection number
of $\delta$ with the medirian $\mu$ agrees with the algebraic
intersection number of $\lambda$ with $\mu$.

In turn, this data gives a map $$\relspinc_{m}\colon \Ta\cap\Tb\longrightarrow
\RelSpinC(Y,K),$$
where $\RelSpinC(Y,K)=\SpinC(Y_0(K))$. To define this map,  
we replace the meridian with a longitude $\lambda$ for the knot,
chosen to wind once along the meridian so that each intersection point $\x$ has a canonical
pair of  closest points $\x'$ and $\x''$. 
Let $(\Sigma,\alphas,\gammas,w)$ denote the
corresponding Heegaard diagram for $Y_0(K)$
(i.e. $\gammas=\lambda\cup\betas_0)$).  We then define
$\relspinc_{m}(\x)$ to be the $\SpinC$ structure in $\SpinC(Y_0(K))$
corresponding
to $\x'$ (or $\x''$).

Now, if $F$ is a Seifert surface for $K$ in $Y$, then it can be closed
off in $Y_0(K)$ to obtain a closed surface ${\widehat F}$. We can
express the evaluation of $c_1(\spinc_w(\x'))$ on ${\widehat F}$ in
terms of data on the Heegaard diagram. Specifically,
recall~\cite{HolDisk} that each two-dimensional homology class for
$Y_0$ has a corresponding ``periodic domain'' $P$ in $\Sigma$, i.e. a
chain in $\Sigma$ whose local multiplicity at $w$ is zero, and which
bounds curves among the $\alphas$ and $\gammas$. Let $\y\in\Ta\cap\Tc$
be an intersection point, let ${\overline n}_{\y}(P)$ denote the sum
of the local multiplicities of $P$ at the points $y_i$ comprising the
$g$-tuple $y_i$ (taken with a suitable fraction if it lies on the
boundary of $P$, as defined in Section~\ref{HolDiskTwo:sec:Adjunction}
of~\cite{HolDiskTwo}), and let $\chi(P)$ denote the ``Euler measure''
of $P$ (which is simply the Euler characteristic of $P$, if all its
local multiplicities are zero or one). Then, it is shown in
Proposition~\ref{HolDiskTwo:prop:MeasureCalc} of~\cite{HolDiskTwo}
that
\begin{equation}
\label{eq:MeasureCalc}
\langle c_1(\spinc_w(\y)),[{\widehat F}]\rangle = 
\chi(P)+2{\Mass}_{\y}(P),
\end{equation}
provided that $P$ is the periodic domain representing ${\widehat F}$.

Returning to the Heegaard diagram for the knot $(Y,K)$, there is a
complex $\CFKinf(Y,K,F)$ generated by tuples
$[\x,i,j]\in\left(\Ta\cap \Tb\right)\times \Z\times \Z$ satisfying the
constraint that 
$$\langle c_1(\relspinc_\Mark(\x)),[{\widehat F}]\rangle + 2(i-j)=0,$$ and with
boundary operator given by $$\partial[\x,i,j]=\sum_{\y}
\sum_{\phi\in\pi_2(\x,\y)}\#\left(\frac{\ModFlow(\phi)}{\R}\right)
[\y,i-n_w(\phi),j-n_z(\phi)]$$
and $\Z[U]$ action defined by
$$U\cdot[\x,i,j]=[\x,i-1,j-1].$$
This complex has a $\Z^2$-filtration given by
$$\Filt[\x,i,j]=(i,j).$$
In particular, $\CFKinf(Y,K,F)$ is generated as a $\Z[U,U^{-1}]$-module
by intersection points $\Ta\cap\Tb$.

As an example, let $T_r$ denote the right-handed trefoil.  The chain
complex $\CFKinf(S^3,F,T_r)$ is generated (as a $\Z[U,U^{-1}]$ module) by
three generators $x$, $y$, and $z$. The filtration degrees are given by
$$\Filt(x)=(-1,0),\hskip.2in \Filt(y)=(0,-1),
\hskip.2in \Filt(z)=(0,0),$$ and the boundary operator is
given by $\partial z=x+y$, and $\partial x=\partial y=0$.

Dually, for the left-handed trefoil $T_{\ell}$ with the torsion $\SpinC$ 
structure, the chain complex is generated by $x$, $y$, $z$ with 
$$\Filt(x)=(1,0),\hskip.2in \Filt(y)=(0,1),
\hskip.2in \Filt(z)=(0,0),$$ and the boundary operator is
given by $\partial x = \partial y = z$, and $\partial z = 0$.  (For
more, see Proposition~\ref{Knots:prop:Trefoil} of~\cite{Knots}.)

In our present applications we will be working with $\HFa(Y)$ (rather
than the full $\HFinf(Y)$), so we can avoid using the full knot
complex.  Rather, consider the naturally induced complex
$\CFK^{0,*}(Y,K)$ generated by $[\x,0,j]$ with $$\langle
c_1(\relspinc_\Mark(\x)),[{\widehat F}]\rangle -2j=0.$$ As a chain
complex, this is the complex for $\CFa(Y)$ using the reference point
$w$. The extra reference point $z$ endows this complex with a relative
$\Z$ filtration; the data of the Seifert surface $F$ comes into play
when one lifts this to an absolute $\Z$-filtration.  For a given
integer $j$, we let $\FiltZero(Y,K,F,j)$ be the subcomplex of
$\CFK^{0,*}(Y,K,F)\cong
\CFa(Y,K,F)$ with filtration degree $\leq j$.

There is an even simpler knot invariant, 
$\CFKa(Y,K,j)$, which is generated by intersection points $\x$
with
$$\langle c_1(\relspinc_{\Mark}(\x)),[{\widehat F}]\rangle = 2j,$$
and with boundary maps counting only those
holomorphic disks with $n_w(\phi)=n_z(\phi)=0$. This complex can,
alternately, be viewed as an associated graded complex for 
$\CFK^{0,*}(Y,K,F)$, associated to its filtration.

\subsection{Notational remarks}
\label{subsec:Notation}

The filtered complexes $\CFKinf(Y,K,F)$ and $\CFK^{0,*}(Y,K,F)$
described here can be written in the notation of~\cite{Knots}, as
follows:
\begin{eqnarray*}
\CFKinf(Y,K,F)&=& \bigoplus_{\{\relspinct\in\RelSpinC(Y,K)\big| \langle
c_1(\relspinct),[{\widehat F}]\rangle = 0\}}\CFKinf(Y,K,\relspinct) \\
\CFK^{0,*}(Y,K,F)&=& \bigoplus_{\{\relspinct\in\RelSpinC(Y,K)\big| \langle
c_1(\relspinct),[{\widehat F}]\rangle = 0\}}\CFK^{0,*}(Y,K,\relspinct). 
\end{eqnarray*}

When $F$ is understood from the context or irrelevant (e.g. when $Y$ is a
rational homology three-sphere), we drop it from the notation.

\subsection{Knot homologies and connected sums}
\label{subsec:ConnSums}

In~\cite{HolDisk}, we give a multiplication map
\begin{equation}
\label{eq:ConnSumMap}
\CFa(Y_1)\otimes_\Z \CFa(Y_2)\longrightarrow \CFa(Y_1\# Y_2),
\end{equation}
which induces a homotopy equivalence. 

In~\cite{Knots}, we give a refinement for the knot invariant, a
variant of which we state in the next proposition. To set up notation,
suppose that $(Y_1,K_1,F_1)$ and $(Y_2,K_2,F_2)$ is a pair of three-manifolds
equipped with oriented null-homologous knots and corresponding Seifert surfaces.
We can form the connected sum
$(Y_1\# Y_2,K_1\# K_2,F_1\#_b F_2)$, where the Seifert surface $F_1\#_b F_2$ is
obtained by a boundary connected sum of the Seifert surfaces $F_1$ and $F_2$. 

\begin{prop}
\label{prop:ConnSumNaturality}
Let $(Y_1,K_1)$ and $(Y_2,K_2)$ be a pair of oriented knots.  
The connected sum map mentioned above induces 
a homotopy equivalence of $\Z$-filtered complexes 
$$
\CFK^{0,*}(Y_1,K_1,F_1)\otimes_\Z
\CFK^{0,*}(Y_2,K_2,F_2)\longrightarrow \CFK^{0,*}(Y_1\# Y_2,K_1\#
K_2,F_1\#_b F_2). $$
\end{prop}

\begin{proof}
This follows readily from Theorem~\ref{Knots:thm:ConnectedSumsOfKnots}
of~\cite{Knots}.
\end{proof}

\begin{remark}
In the statement of Proposition~\ref{prop:ConnSumNaturality}, the $\Z$
filtration on the left-hand-side comes from the following convention:
the tensor product of two $\Z$-filtered complexes $(C_1,\Filt_1)$ and
$(C_2,\Filt_2)$ is naturally a $\Z$-filtered complex $(C_1\otimes
C_2,\Filt_{\otimes})$, where the tensor product filtration is defined
by $$\Filt_{\otimes}(x_1\otimes x_2)=\Filt_1(x_1)+\Filt_2(x_2)$$ for
arbitrary homogeneous generators $x_1\in C_1$ and $x_2 \in C_2$.
\end{remark}

\section{Invariants of fibered knots}
\label{sec:FiberedKnots}

We now prove Theorem~\ref{intro:FiberedKnots}. The alternate
description of $c(K)$ which involves only Floer homology groups for
closed three-manifolds which fiber over the circle (rather than the
knot invariant) will follow quickly from the proof.

\vskip.2cm
\noindent{\bf{Proof of Theorem~\ref{intro:FiberedKnots}.}}
The proof proceeds by constructing, for any fibered knot $(Y,K)$,
a genus $2g+1$ Heegaard diagram with the following properties:
\begin{list}
        {(H-\arabic{bean})}{\usecounter{bean}\setlength{\rightmargin}{\leftmargin}}
\item 
\label{item:Uniqueness}
there is a unique intersection point $\x\in\Ta\cap\Tb$ with
$$\langle c_1(\relspinc_\Mark(\x)), [{\widehat F}]\rangle = -2g,$$
\item 
\label{item:Adj}
there are no intersection points $\y\in\Ta\cap\Tb$ with
$$\langle c_1(\relspinc_\Mark(\y)), [{\widehat F}]\rangle < -2g,$$
\item
\label{item:Admissibility}
the Heegaard diagram for the knot is is weakly admissible
for any $\SpinC$ structure over $Y$.
\end{list}
(Recall that a pointed Heegaard diagram is weakly admissible for any
$\SpinC$ structure if all the periodic domains have both positive and
negative local multiplicities. For the definition and discussion of
this hypothesis, see~\cite{HolDisk}.) Given such a diagram, the
theorem follows at once.

The Heegaard diagram is constructed as follows. Start with a genus $g$
surface $A$ with boundary consisting of two circles, denoted
$\alpha_{1}$ and $\lambda$. We can find two $2g$-tuples of pairwise
disjoint arcs in $A$, $(\xi_2,...,\xi_{2g+1})$ and
$(\eta_2,...,\eta_{2g+1})$, and an arc $\delta$ connecting $\alpha_1$
to $\lambda$ with
the following properties:
\begin{list}
        {(\arabic{bean})}{\usecounter{bean}\setlength{\rightmargin}{\leftmargin}}
\item 
\label{item:Disjointness}
The curve $\xi_i$ is disjoint from $\eta_j$ unless $i=j$,
in which case $\xi_i$ meets $\eta_i$ transversally in a single intersection point
which we denote by $x_i$.
\item 
The curve $\delta$ is disjoint from all the $\xi_i$ and $\eta_i$.
\item 
If we let ${\overline A}$ denote a copy of $A$ reflected across its boundary,
and let $\Sigma$ denote the closed surface of genus $2g+1$ obtained by
gluing $A$ to ${\overline A}$ along their boundary, let $\alpha_i$
(for $i=2,...,2g+1$) denote the closed curve obtained by gluing $\xi_i$
to ${\overline \xi}_i$, and similarly let $\beta_i$ denote the closed
curve obtained by gluing $\eta_i$ to ${\overline \eta}_i$, then
$(\Sigma,\{\alpha_1,...,\alpha_{2g+1}\},\{\lambda,\beta_2,...,\beta_{2g+1}\})$
is a Heegaard diagram for $S^1\times {\widehat F}$.
\item Indeed, letting $\mu$ denote the closed curve obtained by  joining $\delta$ to
${\overline \delta}$, the diagram
$(\Sigma,\alphas,\{\beta_2,...,\beta_{2g+1}\},\mu)$ represents the
fibered knot $B\subset \#^{2g}(S^2\times S^1)$ whose monodromy map is
the identity map.
\item There are pairwise disjoint curves $\sigma_1,...,\sigma_g$ 
with the property that $\sigma_i\cap \xi_j=\emptyset$ 
except when $j=2i+1$, in which case $\sigma_i$ 
meets $\xi_{2i+1}$ in a single transverse intersection point.
\item There are pairwise disjoint curves $\tau_1,...,\tau_g$ with the property that
$\tau_i\cap\xi_j=\emptyset$ except when $j=2i$, in which case
$\tau_i$ meets $\xi_{2i}$ in a single, transverse double-point. Moreover,
$\tau_i\cap\eta_j=\emptyset$, except when $j=2i+1$, in which
case $\tau_i$ meets $\eta_j$ in a single, transverse double-point.
\end{list}

\begin{figure}
\mbox{\vbox{\epsfbox{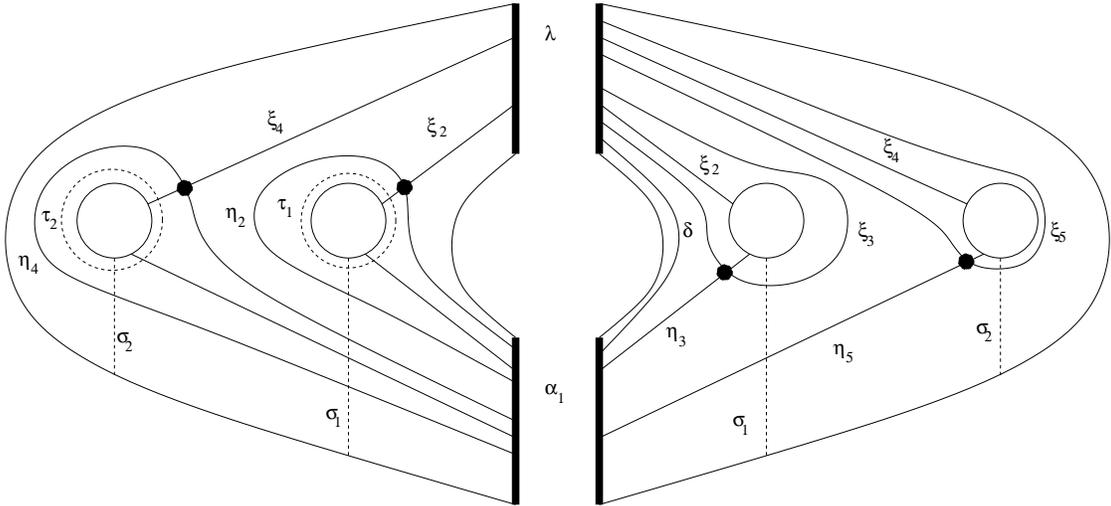}}}
\caption{\label{fig:RegionA}
{\bf{Region A.}}  An illustration of the region $A$ described in the
proof of Theorem~\ref{intro:FiberedKnots}, in the case where
$g=2$. Note that we have drawn the ``top'' and the ``bottom'' of $A$:
the region $A$ is obtained by identifying the two planar regions
illustrated along the two circles and two light arcs (the dark arcs
labeled $\alpha_1$ and $\lambda$ are not identified in in their
interiors; rather, they glue up to make closed curves in $A$).  The
intersection points $x_i$ for $i=2,...,5$ are marked with heavy dots
(but not labelled).}
\end{figure}

In particular, we have a Heegaard diagram for $S^1\times{\widehat F}$ which is
divided in two by $\alpha_{1}\cup\lambda$.  Note that the
Heegaard splitting for $S^1\times {\widehat F}$ is the one obtained by
first dividing the circle into two intervals $I_1$ and $I_2$,
attaching a one-handle to $I_1\times {\widehat F}$, which is drilled
out of $I_2\times {\widehat F}$, and then symmetrically attaching a
one-handle to what is left of $I_2\times {\widehat F}$ which in turn is
drilled out of $I_1\times {\widehat F}$. For an illustration of the region $A$,
see Figure~\ref{fig:RegionA}

\begin{figure}
\mbox{\vbox{\epsfbox{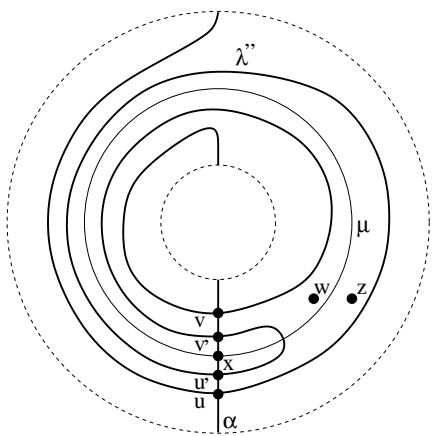}}}
\caption{\label{fig:WindTwice}
{\bf{Winding about $\mu$.}}  An illustration of the modifications to
the Heegaard diagram made in an annular neighborhood of
$\mu=\delta\cup{\overline \delta}$. 
We have dropped some subscripts from the picture:
in the picture, $x$ and $\alpha$ denote $x_1$ and $\alpha_1$ from
the text.}
\end{figure}

After a suitable modification of the Heegaard, diagram, we would like
to place a basepoint on the meridian $\mu$. To this end, wind
$\lambda$ half way along $\mu$, so that the new curve $\lambda'$ meets
$\alpha_{1}$ in two intersection points $u$ and $v$. Corresponding to
this finger move, there is now a new region $\cald$ in the Heegaard
surface bounded by an arc in $\alpha_{1}$ and an arc in $\lambda'$. We
place the basepoint $m$ on $\mu\cap \cald$. This basepoint gives rise
then to a pair $w$ and $z$ of basepoints in the Heegaard surface
$\Sigma$. Without crossing $w$ or $z$, we wind one more time
completely along $\mu$ to get $\lambda''$, so that the periodic domain
$P$ bounded by $\alpha_1$ and $\lambda''$ has both positive and
negative coefficents, i.e. we introduce a new region $\cald'$ bounded
by an arc in $\alpha_1$ and $\lambda''$, and let $u'$ and $v'$ denote
the newly-introduced pair of intersection points which lie on the
boundary of the closure of $\cald'$.  Let $x_{1}$ be the intersection
point between $\mu$ and $\alpha_{1}$ (see Figure~\ref{fig:WindTwice}
for an illustration). Clearly, given any $\y\in\Ta\cap\Tb$, we have
that $x_1\in\y$; and indeed, if we replace $x_1\in\y$ by $u'$, we
obtain an intersection point $\y'\in \Ta\cap\Tc$ (where here
$\gammas=\{\lambda'',\beta_2,...,\beta_{2g+1}\}$) which represents
$\relspinc_\Mark(\y)$.

Let $P$ denote the periodic domain with $\partial
P=\alpha_1+\lambda$.  Clearly, $P$ represents the homology
class of ${\widehat F}$. By construction, its local multiplicities are
$-1$ (in $\cald'$), $0$ in the winding region, $1$ in the region
corresponding to $A$ and $2$ (in the region corresponding to
${\overline A}$). 
It is now an easy consequence of
Property~\eqref{item:Disjointness} that:
\begin{list}
        {(\arabic{bean}')}{\usecounter{bean}\setlength{\rightmargin}{\leftmargin}}
\item 
\label{item:UniqueMaximal}
If we let $\x=x_1\times...\times x_{2g+1}$, 
then for any other point $\y\in\Ta\cap\Tb$ (i.e. with $\y\neq \x$),
$$2g=\Mass_{\x'}(P)<\Mass_{\y'}(P).$$
\end{list}
As a Heegaard diagram for $\#^{2g}(\#S^2\times S^1)$,
$(\Sigma,\alphas,\betas'',w)$ is not weakly admissible yet: periodic
domains representing the $2g$ tori have only non-negative local
multiplicities. However, we can arrange for the diagram to be
admissible by winding along the $2g$ curves
$\tau_1,...,\tau_g$ and 
${\overline \sigma}_1,...,{\overline \sigma}_g$ (of course, here
$\tau_i$ are curves in $A\subset A\cup {\overline A}=\Sigma$, while
${\overline \sigma}_i$ are curves in ${\overline A}\subset \Sigma$
corresponding to the $\sigma_i$ described above). 

We claim that this procedure does
not affect Property~(\ref{item:UniqueMaximal}').
Specifically, since all the  curves
$\{\tau_i,{\overline \sigma}_i\}_{i=1}^g$
are disjoint from $\alpha_1$ and
$\lambda$, it follows that the periodic domain $P$ is unaffected by
this winding.  Thus, it remains the case that if $\y\in\Ta\cap\Tb$, then
$2g=\Mass_{\x'}(P)\leq \Mass_{\y'}(P)$.  Verifiying that
Property~(\ref{item:UniqueMaximal}') still holds amounts to showing
that if $\Mass_{\x'}(P)=\Mass_{\y'}(P)$ then $\x=\y$.  To see that this
holds after winding along the $\tau_i$, we let $\y\in\Ta\cap\Tb$
satisfy $2g=\Mass_{\y'}(P)$. Then, since the $\tau_i$ are disjoint from
all $\xi_{2j}$ except when $i=j$, we see that $x_{2i}\in\y$, while the
fact that the only $\xi$-curve met by $\tau_i$ is $\xi_{2i}$ while the
only $\eta$-curve met by $\tau_i$ is $\eta_{2i+1}$ ensures that none
of the newly-introduced intersection points between $\xi_i'$ and
$\eta_j$ appears in the tuple $\y$.  Thus, $\y=\x$. The fact that this
property is preserved after winding along the
$\{{\overline\sigma}_i\}$ follows easily from the fact that those
curves are supported in ${\overline A}$, a region where the local
multiplicty of $P$ is $2$ (instead of $1$).

In this Heegaard diagram for $(Y,K)=(\#^{2g}(S^2\times S^1),B)$, we
have achieved admissibility
(Property~(H-\ref{item:Admissibility})). The other two required
properties -- Property~(H-\ref{item:Uniqueness}) and
(H-\ref{item:Adj}) now follow readily from
Property~(\ref{item:UniqueMaximal}'), together with
Equation~\eqref{eq:MeasureCalc} (noting that the ``Euler measure'' of
the periodic domain $P$ is $-6g$), so that $$\langle
c_1(\relspinc_\Mark(\x),[{\widehat F}]\rangle=\langle
c_1(\spinc_{w}(\x')),[{\widehat F}] \rangle = -2g.  $$ 
(Where here again the point
$\x'$ is obtained by
substiting $u'$ for the $x_{1}$ coordinate in $\x$.)

In general, suppose we have an open book decomposition with monodromy
map $\phi$.  We can extend $\phi$ to an automorphism $\Phi$ of
$\Sigma$, by viewing the complement of a neighborhood of ${\overline
\delta}$ in ${\overline A}$ as $F$, and extending $\phi$ over the rest
$\Sigma$ as the identity map. It is easy to see that the Heegaard
diagram $(\Sigma,\alphas,\mu\cup \Phi(\beta_2,...,\beta_{2g+1}))$
represents the given open book. Weak admissibility is still easy to
see, and indeed Property~(\ref{item:UniqueMaximal}') still holds (the
argument in the model case still applies, since the general Heegaard
diagram and the model case agree in the region of $\Sigma$ where the
multiplicities of the periodic domain $P$ are less than two).
\qed
\vskip.3cm

Suppose $Y_0$ is a three-manifold with a distinguished two-dimensional
homology class $[{\widehat F}]$. In this case, we define
$$\HFp(Y_0,{\widehat F},i)
=
\bigoplus_{\{\spinc\in\SpinC(Y_0)\big|
\langle c_1(\spinc),[{\widehat F}]\rangle = 2i\}}
\HFp(Y_0,\spinc).
$$
It is shown in~\cite{HolDiskSymp}, 
that if $Y_0$ is a three-manifold whose
fiber ${\widehat F}$ has genus $g>1$,
then 
\begin{equation}
\label{eq:ThreeManifoldsFiber}
\HFp(Y_0,{\widehat F},1-g)
\cong \Z.
\end{equation}

\begin{prop}
\label{prop:FiberedCharacterization}
Let $(Y,K)$ be a fibered knot with genus $g>1$ and let $F$ be a
Seifert surface. Let $c(K)$ denote the image of a generator of
$H_*(\FiltZero(-Y,K,F,-g))$ inside $\HFp(-Y)$.  Moreover, let
${\widehat c}(K)$ be a generator of $\HFa(-Y_0(K),{\widehat F},1-g)$
whose image in $\HFp(Y_0,{\widehat F},1-g)$ is a generator, and let
$$F_W\colon \HFa(-Y_0)\longrightarrow \HFa(-Y)$$ denote the map
induced by the cobordism (where we suppress the $\SpinC$ structure, as
it is determined by its restriction to $-Y_0$). Then $$F_{W}({\widehat
c}(K))=\pm c(K).$$
\end{prop}

\begin{figure}
\mbox{\vbox{\epsfbox{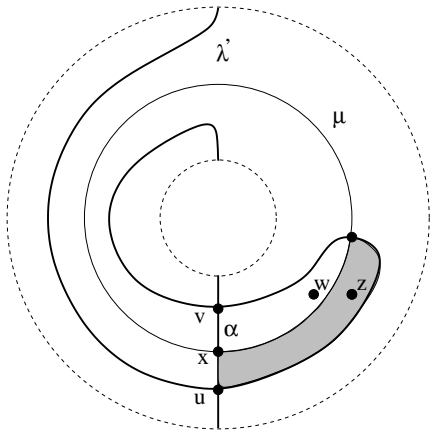}}}
\caption{\label{fig:WindOnlyOnce}
{\bf{Illustration of the proof of
Proposition~\ref{prop:FiberedCharacterization}.}}  Again, an
illustration of a neighborhood of $\mu$. This time, we wind
$\lambda$ only once, however, to obtain the curve $\lambda'$
(denoted $\beta'$ here). The domain of small triangle $\psi$ representing
the map from the canonical generator $-Y_0$ to $-Y$ is
shaded.}
\end{figure}

\begin{proof}
Note that in the proof of Theorem~\ref{intro:FiberedKnots}, the
Heegaard diagram for $Y_0(K)$ is weakly admissible for any
$\SpinC$ structure. Indeed, consider the Heegaard diagram where we
undo the last finger move which created $\cald'$ (thus cancelling the
two closest points $u'$ and $v'$ -- note that we used a diagram with
this extra pair of intersection points in the previous proof only to
calculate $\langle c_1(\relspinc_\Mark(\x)),[{\widehat F}]\rangle$;
however, now that we know this is $-2g$, we no longer need the extra
pair of intersection points). This diagram is still weakly admissible
for any $\SpinC$ structure whose evaluation on ${\widehat F}$ is
non-trivial.  Moreover, for this Heegaard diagram, it is easy to see
(using Equation~\eqref{eq:MeasureCalc}) that there are exactly two
intersection points ${\mathbf u}=u\times x_2\times...\times x_{2g+1}$
and ${\mathbf v}=v\times x_2\times...\times x_{2g+1}$ which represent
a $\SpinC$ structure $\spinc$ satisfying $$\langle c_1(\spinc),
[{\widehat F}]\rangle = 2-2g.$$ In fact, using a domain $\cald$
supported in the winding region, we see that there is a unique
$\phi\in\pi_2({\mathbf u},{\mathbf v})$ with $\Mas(\phi)=1$,
$n_w(\phi)=1$, and $\#\ModFlow(\phi)=1$. Indeed, this gives a direct
proof of Equation~\eqref{eq:ThreeManifoldsFiber}\footnote{Recall that
Equation~\eqref{eq:ThreeManifoldsFiber}, which can be thought of the
analogue of Theorem~\ref{intro:FiberedKnots} for fibered
three-manifolds, is established in~\cite{HolDiskSymp} using general
properties of $\HFp$: adjunction inequalities, surgery long exact
sequences, and one model calculation. In fact,
Theorem~\ref{intro:FiberedKnots} can be proved this way as well, 
using corresponding properties of the knot invariants.}, giving
also a geometric representative $[{\mathbf u},0]$ for the cycle
representing the generator of $\HFp(Y_0,{\widehat F},1-g)$ (and hence
${\mathbf u}$ represents the corresponding ${\widehat c}(K)\in
\HFa(-Y_0(K),1-g)$).

For the Heegaard diagram, there is also a small triangle
$\psi\in\pi_2({\mathbf u},\Theta,\x)$ with $\Mas(\psi)=0$,
$\#\ModFlow(\psi)=1$, $n_w(\psi)=0$ and $n_z(\psi)=1$.  It follows
that there is no $\psi'\in\pi_2({\mathbf u},\Theta,\y)$ with
$n_w(\psi')=0$, $\Mas(\psi')=0$, and $\cald(\psi')\geq 0$. To see
this, observe that the topological constraints on the cobordism show
that any triangle $\psi'\in\pi_2({\mathbf u},\Theta,\y)$ can be
decomposed as $\psi'=\psi *\phi$, where $\phi\in\pi_2(\x,\y)$. Now,
according to properties of the Heegaard diagram for $(Y,K)$
constructed in the proof of Theorem~\ref{intro:FiberedKnots}
(specifically, Properties~(H-\ref{item:Uniqueness}) and
(H-\ref{item:Adj}), it follows that if $\x\neq \y$, then
$n_z(\phi)=0$, contradicting $\cald(\psi')>0$. In the case where
$\x=\y$, then, the Whitney disk $\phi$ represents a periodic domain,
with $\Mas(\phi)=0$. According to admissibility, then, $\cald(\phi)$
has both positive and negative coefficients, so it follows readily
that $\psi'$ has this same property.

It now follows readily that $F_{W}({\mathbf u})=\x$, as
claimed. (See Figure~\ref{fig:WindOnlyOnce} for an illustration.)
\end{proof}


We conclude with a few examples.

The unknot $U\subset S^3$ has a genus one Heegaard diagram, with a
single intersection point. It is obvious in this case that
$H_*(\FiltZero(-S^3,U,0))\cong\Z$, and that, $c(U)$ is a
generator for $\HFa(-S^3)\cong \Z$.  We have the following more
interesting example, which will also help in the proofs of
Theorems~\ref{intro:OverTwisted} and
\ref{intro:SteinFillable}.

\begin{lemma}
\label{lemma:Trefoil}
Let $T_\ell$ resp. $T_r$ denote the left- resp. right-handed trefoil
inside $S^3$. 
Then, $c(T_r)\in \HFa(-S^3)\cong \Z$ 
is a generator, 
whereas 
$c(T_\ell)=0$.
\end{lemma}

\begin{proof}
This is an immediate consequence of the calculation of $\CFKinf$ for the
trefoil given in Proposition~\ref{Knots:prop:Trefoil} of~\cite{Knots},
see also Subsection~\ref{subsec:KnotHomology} above.
\end{proof}

\begin{remark}
A generalization of Lemma~\ref{lemma:Trefoil} to alternating, fibered knots is
given in~\cite{AltKnots} (see especially 
Corollary~\ref{AltKnots:cor:AltContact} of~\cite{AltKnots}).
\end{remark}

\section{Invariants of contact structures}
\label{sec:Equivalence}

We now prove Theorem~\ref{intro:SameContact}: if two fibered knots
induce the same contact structure, then their invariants coincide. The
proof is broken into several lemmas. After establishing this theorem,
we explain (in Subsection~\ref{subsec:HomotopyClass}) the relationship
between $c(Y,K)$ and the homotopy class of the two-plane field
underlying $\xi$.

\begin{lemma}
\label{lemma:IndepTref}
For the Giroux stabilization of a fibered knot $(Y,K)$
obtained by attaching a trefoil to $K$, we have that
$c(Y,K)=c(Y,K\#T_r)$.
\end{lemma}

\begin{proof}
This follows
from the naturality under connected sums, together with the sample
calculation from Lemma~\ref{lemma:Trefoil}. More precisely, by
Proposition~\ref{prop:ConnSumNaturality}, we see that 
$$c(K\# T_r)=
\iota(c_0(K\# T_r))=\iota(c_0(K)\otimes c_0(T_{r}))= 
c(K)\otimes c(T_{r}),$$
where the last tensor product indicates the isomorphism
$$\HFa(-Y)\otimes \HFa(-S^3)\longrightarrow \HFa(-Y).$$  But by
Lemma~\ref{lemma:Trefoil}, $c(T_r)$ generates $\HFa(-S^3)$, so it follows
that $c(K\# T_r)=c(K)$.
\end{proof}

We establish now a naturality for $c(K)$ under Dehn twists, which will
be used repeatedly.

\begin{theorem}
\label{thm:DehnNaturality}
Let $(M,\phi)$ be an open book decomposition for $Y$, and fix a curve
$\gamma\subset Y-K$ supported in a page of the open book, which is not
homotopic to the boundary. Then, 
$(M,\phi\circ\RightDehn_\gamma^{-1})$ induces an
open book decomposition of $Y_{+1}(\gamma)$, and indeed under the map
$$F_{W}\colon \HFa(-Y)\longrightarrow \HFa(-Y_{+1}(\gamma))$$
obtained by the two-handle addition (and summing over all $\SpinC$ structures),
we have that
$$F_{W}(c(M,\phi))=\pm c(M,\phi\circ\RightDehn_\gamma^{-1}).$$
\end{theorem}

\begin{proof}
In view of Lemma~\ref{lemma:IndepTref}, it suffices to consider the
case where the genus $g$ of the fiber is $>1$.  Let $Y_0=Y_0(K)$,
$Y'=Y_{+1}(\gamma)$, and let $Y'_0$ denote the three-manifold obtained
by performing a zero-surgery on $Y'$ along the knot $K'\subset Y'$
induced from $K$.  Also, let $V$ and $V'$ denote the natural
cobordisms from $-Y_0$ to $-Y$ and $-Y_0'$ to $-Y'$ respectively and let
$W_0$ be the cobordism from $-Y_0$ to $-Y_0'$ obtained by attaching the
two-handle along $\gamma$, thought now as a curve in $Y_0$ (rather
than $Y$). Then, by naturality of the maps induced by cobordisms
(see~\cite{HolDiskFour}), we have a commutative diagram $$
\begin{CD}
\Z\cong \HFa(-Y_0,1-g) @>{F_{W_0}}>> \HFa(-Y'_0,1-g) \cong \Z\\
@V{F_{V}}VV @VV{F_{V'}}V \\
\HFa(-Y) @>{F_W}>> \HFa(-Y').
\end{CD}
$$ Now, $F_{W_0}$ induces an isomorphism, since the map fits into an
exact sequence where the third term is a three-manifold for which the
corresponding Floer homology groups vanish by the adjunction
inequality (c.f. Section~\ref{HolDiskSymp:sec:NonVanishing}
of~\cite{HolDiskSymp}).  The theorem now follows immediately from
Proposition~\ref{prop:FiberedCharacterization}.
\end{proof}

\begin{remark}
Note that the use of Proposition~\ref{prop:FiberedCharacterization} after
stabilizing can be avoided (both here and in the proof of
Theorem~\ref{intro:SteinFillable} below) by using corresponding
properties of the knot invariant directly.
\end{remark}

\begin{lemma}
\label{lemma:DependsOnlyOnGenus}
Let $(Y,K)$ be a fibered knot specified by an open book $(F_g,\phi)$. Then, for each
integer $h\geq 0$, there
is a class $c(\phi,h)\in \HFa(-Y)$ with the property that for any 
genus $h$ Giroux stabilization $\phi'$ of $(F_g,\phi)$, $c(\phi')=c(\phi,h)$.
\end{lemma}

\begin{proof}
Note that the $\phi\#_b\Id_h\in \MClass(F_g\#_b F_h)$ induces an open book
decomposition of 
$Y\# \left(\#^{2h}(S^2\times S^1)\right)$. By attaching $2h$ three-handles to 
this latter three-manifold, we obtain a cobordism $U$ to $Y$. Let $c(\phi,h)$ be the class
obtained as $F_{U}(c(\phi\#_b\Id_h))$.

Now, let $\phi'$ be any Giroux stabilization of $\phi$. Then, we can
``untwist'' $\phi'$ by composing it with $2h$ 
left-handed Dehn twists
to obtain $\phi\#_b \Id_h$ (i.e. we take here left-handed Dehn twists
along parallel copies of the curves used for the Giroux
stabilization). By Theorem~\ref{thm:DehnNaturality}, the image of
$c(\phi')$ under the induced $2h$ two-handles $W$ gives $c(\phi\#_b
\Id)$. We claim that the three-handles in $U$ cancel the two-handles
of $W$, so that $U\circ W\cong [0,1]\times Y$, and hence that
$$c(\phi')=F_U(c(\phi\#_b\Id))$$ is independent of the particular
stabilization $\phi'$.
\end{proof}

\vskip.2cm
\noindent{\bf{Proof of Theorem~\ref{intro:SameContact}.}}
According to Giroux's theorem, we need to verify only that the class
$c(K)$ is invariant under Giroux stabilizations. But this follows
immediately from Lemmas~\ref{lemma:IndepTref} and~\ref{lemma:DependsOnlyOnGenus}. 
\qed

\begin{remark}
Observe that when $b_1(Y)>0$, it is an easy consequence of
Theorem~\ref{intro:FiberedKnots} that the invariant $c(K)$ lies in the
kernel of the action of $H_1(Y;\Z)$. In a related direction, the
constructions descibed here carry over verbatim to the case of twisted
coefficients, giving rise to an element ${\underline
c}(K)\in\uHFa(-Y)$, whose image under the natural map
$\uHFa(-Y)\longrightarrow \HFa(-Y)$ is the element $c(K)$ described
above. We have no further use for this refinement in the present
paper.
\end{remark}

\subsection{Classical data}
\label{subsec:HomotopyClass}

Recall that an oriented two-plane field $\xi$ in an oriented
three-manifold $Y$ has an induced $\SpinC$ structure, the equivalence
class of the positively oriented normal vector, which we denote
$\spinc(\xi)$. When the first Chern class of this $\SpinC$ structure
is torsion, the homotopy type of $\xi$ is uniquely determined by the
pair $\spinc(\xi)$, and an additional Hopf-type invariant $h(\xi)\in \Q$,
which we describe presently (see~\cite{GompfStipsicz} and
\cite{KMcontact} for a detailed discussion). Suppose that $(W,J)$ is an
almost-complex four-manifold which bounds $Y$ so that the complex
tangencies at its boundary agree with the two-plane field $\xi$. Then,
the Hopf invariant (which is independent of the $\SpinC$ structure) is
given by $$h(\xi)=\frac{c_1(W,J)^2+2-2\chi(W)-2\sigma(W)}{4},$$ where
here $c_1(W,J)^2$ denotes the evaluation on $[W,\partial W]$ of the
square of the first Chern class of the almost-complex structure,
$\chi(W)$ denotes the Euler characteristic of $W$, and $\sigma(W)$
denotes the signature of its intersection form.

\begin{prop}
The class $c(\xi)\in\HFa(Y)$ is supported in the summand
$\HFa(Y,\spinc(\xi))$ corresponding to the $\SpinC$ structure
belonging to the contact structure.  Moreover, if $c_1(\spinc(\xi))$
is a torsion class, then the absolute grading of $c(\xi)$ agrees with
the Hopf invariant $h(\xi)$.
\end{prop}

\begin{proof}
We assume without loss of generality that the genus $g$ of the open
book decomposition representing $\xi$ is greater than one. (Note that
stabilization fixes both the $\SpinC$ structure of $c(\xi)$ and its
absolute grading).  Now, in~\cite{HolDiskSymp} it is shown that the
generator for $\HFp(Y_0,1-g)$ is supported in the canonical $\SpinC$
structure of the fibered three-manifold.  From this and
Proposition~\ref{prop:FiberedCharacterization}, the result on $\SpinC$
structures follows immediately.

For the statement about absolute gradings, we use the fact that the
two-handle addition from $-Y_0$ to $-Y$ can be given an almost-complex
structure (this follows from a local analysis). In~\cite{HolDiskSymp},
we constructed, for each three-manifold which fibers over the circle,
a Lefschetz fibration $W_0$ over the disk, which we think of now as a
cobordism from $S^3$ to $-Y_0$, with the property that, if we endow
$W$ with the $\SpinC$ structure associated to its canonical
almost-complex structure, then $F_{W_0,\spinccan}(\theta)\in
\HFp(Y_0,1-g)\cong \Z$ is a  generator (where here $\theta\in\HFp(S^3)$ is 
a zero-dimensional generator).  Following
Proposition~\ref{prop:FiberedCharacterization}, $c(\xi)$ is the image
of $F_{W_0\cup W}$, endowed with its canonical $\SpinC$ structure. The
absolute grading of this element is immediately seen to be given by
$h(\xi)$.
\end{proof}

\section{Vanishing and non-vanishing results}

The vanishing theorem for overtwisted contact structures, and the
non-vanishing theorem for Stein fillable ones, now follow rather
easily from what has already been established.

\vskip.2cm
\noindent{\bf{Proof of Theorem~\ref{intro:OverTwisted}.}}
If $(Y,K)$ is an arbitrary fibered knot, we prove that $(Y,K\#
T_\ell)$ has vanishing invariant.  This follows immediately from the
naturality of the knot invariant under connected sums
(Proposition~\ref{prop:ConnSumNaturality}), together with the sample
calculation of Lemma~\ref{lemma:Trefoil}. 
This suffices to establish
the result, in light of Eliashberg's classification of overtwisted 
contact structures~\cite{EliashbergOvertwisted}.
\qed

\vskip.2cm
\noindent{\bf{Proof of Theorem~\ref{intro:SteinFillable}.}}
Let $\phi$ be a word in the mapping class group of $F_g$ corresponding
to a Stein fillable contact structure. By Giroux's characterization,
we can ``untwist'' $\phi$ using only left-handed Dehn twists to
obtain the identity map, so by the naturality from
Theorem~\ref{thm:DehnNaturality}, we get a cobordism $W$ from $-Y$ to
$-\#^{2g}(S^2\times S^1)$ with the property that $$F_{W}(c(\xi))=\pm
c(\xi_0),$$ where $\xi_0$ denotes the contact structure over
$\#^{2g}(S^2\times S^1)$ corresponding to the identity map on
$F_{2g}$. The fact that $c(\xi_0)$ is non-trivial follows as in the
proof of Lemma~\ref{lemma:DependsOnlyOnGenus} from the fact that the
unknot in $S^3$ has non-trivial invariant.
\qed

\section{Tightness of fibered surgeries}
\label{sec:FourDInterp}

We can now prove the following general result 
(compare Theorem~\ref{intro:MinusOneFibered}):

\begin{theorem}
\label{thm:FiberedGen}
Let $K\subset Y$ be a fibered knot with genus $g>1$ in any three-manifold
with $\HFpRed(Y)=0$.  Then, the induced contact structure on $Y_{-1}(K)$
is tight.
\end{theorem}

\begin{proof}
In view of Proposition~\ref{prop:FiberedCharacterization} and
Theorem~\ref{intro:OverTwisted}, it suffices to prove that the
generator of $\HFp(-Y_0(K),1-g)$ in $\HFp(-Y_{-1}(K))$ under the natural
two-handle addition is non-trivial.  But this follows easily from the
surgery long exact sequence
(c.f. Theorem~\ref{HolDiskTwo:thm:GeneralSurgery}
of~\cite{HolDiskTwo}), which reads $$
\begin{CD}
...@>>> \HFp(-Y)@>>>\HFp(-Y_0(K)) @>>>\HFp(-Y_{-1}(K))@>>>...
\end{CD}
$$ Now, since the $U$-action on $\HFp(-Y_0(K),1-g)$
is trivial, while each element of $\HFp(-Y)$ lies in the image of
the $U$-action (since $\HFpRed(-Y)=0$), it follows that elements in
$\HFp(-Y_0(K),1-g)$ inject into $\HFp(-Y_{-1}(K))$, as
desired.
\end{proof}

\vskip.2cm
\noindent{\bf{Proof of Theorem~\ref{intro:MinusOneFibered}.}}
Since $\HFpRed(S^3)=0$, Theorem~\ref{thm:FiberedGen} handles the cases
where $g>1$.  When $g=1$, there are only three fibered knots in
$S^3$: the two trefoils and the figure eight knot, where the theorem
can be checked. Indeed, in all three cases, we claim that
$S^3_{-1}(K)$ is Stein fillable. To this end, fix any one of these
knots $K$, let $\phi$ denote its monodromy, acting on a surface $F$ of
genus one (and fixing its boundary), and let $c$ denote a curve in the
interior of $F$ which is isotopic to the boundary.
Our aim then is to show that 
the mapping class $\phi\cm R_{c}$ can be written as a product of
right-handed Dehn twists.

We consider the case of the left-handed trefoil. In this case, the
mapping class $\phi$ can be written as a product of left-handed Dehn
twists.  Explicitly, if $a$ and $b$ are a pair of circles in $F$
which meet transversally in a single point, the monodromy 
$\phi=R_a^{-1}\cm
R_b^{-1}$. The fact that $R_a^{-1}\cm R_b^{-1}
\cm R_c$ can be written as a product of right-handed Dehn twists
follows immediately from the following easily verified
relation in the mapping class group:
$$(R_b\cm R_a)^6=R_c.$$

Similarly, since the monodromy of the 
figure eight knot can be written as 
$R_a\cm R_b^{-1}=R_a^2\cm (R_a^{-1}\cm R_b^{-1})$,
the result follows from
the corresponding fact for the left-handed trefoil.

The case of the right-handed trefoil is straightforward.
\qed

\commentable{
\bibliographystyle{plain}
\bibliography{biblio}
}

\end{document}